\documentclass[oneside,a4paper]{article}
\usepackage{pp}
\usepackage{pictexwd,dcpic}
\usepackage{multirow}
\usepackage[russian,english]{babel}
\usepackage[cp1251]{inputenc}

\def\3{\!3}
\def\4{\!4}

\def\SS{\mathbb{S}}
\def\DD{\mathbb{D}}
\def\Iso{\mathrm{Iso}}

\begin{document}

\title{The telescopic construction;\\
a microsurvey
}
\author{D. Panov\thanks{is a Royal Society University Research
Fellow}
\ and A. Petrunin\thanks{was partially supported by NSF grant DMS 1309340}
}
\date{}
\maketitle

\section{Introduction}

It is well known that there are cocompact isometric properly discontinuous actions on hyperbolic $n$-space.
On the other hand almost any question about possible quotient spaces of such  actions $\Gamma\acts \HH^n$
remains open for large $n$;
see Section \ref{sec:open}.

The question we are trying to answer is whether any given finitely presented group
is the fundamental group of $\HH^n/\Gamma$
for an action $\Gamma\acts \HH^n$ as above.
We started to think about this at the 16-th G\"okova Conference
and later we were able to give a positive answer for some dimensions.
To our satisfaction, the construction turned out to be useful in a couple of unexpected places.

Our original proof is given by an explicit but slightly convoluted construction
which might be hard to follow.
In this note we will only try to convince the reader that the construction has enough freedom to attain the goal.
(The same problem happens if one considers the following question: \emph{is it possible to build a doll house from the
Lego blocks?} ---
obviously yes,
but the easiest way to prove this is to build a house from Lego blocks by your hands.
This type of proof is given in \cite{panov-petrunin},
and here we try to give some evidence for the ``obviously yes'' answer.)

We also overview the results in  \cite{FP1},  \cite{K1}, \cite{panov-petrunin}.

\parbf{Acknowledgments.} We want to thank Joel Fine, Misha Kapovich and Jos\'e Mar\'{\i}a Montesinos-Amilibia for help.

\section{What can be proved}

Let us denote by $\langle\Tor \Delta\rangle$
the subgroup of group $\Delta$
generated by the elements of finite order.

The proof of the following theorem is discussed in Section~\ref{sec:telescope}

\begin{thm}{Theorem}\label{thm:main}
There is a finitely presented group $\Gamma$
such that for any finitely presented group $G$
there is a finite index subgroup $\Gamma'$ in $\Gamma$
such that $G$ is isomorphic to the quotient $\Gamma'/\langle\Tor \Gamma'\rangle$.

Moreover,
\begin{enumerate}[(i)]
\item\label{thm:main:0}  The group $\Gamma$ can be chosen to be hyperbolic.
\item\label{thm:main:2} The  Coxeter group $\Gamma_{12}\acts\HH^3$ of the regular right angled hyperbolic dodecahedron satisfies this property. Moreover
\begin{enumerate}
\item\label{thm:main:2:a} The same holds for a subgroup of finite index in $\Gamma_{12}$ which does not contain reflections in planes.
\item\label{thm:main:2:b}
The same holds for a subgroup infinite index
in $\Gamma_{12}$ which acts cocompactly on a convex subset of $\HH^3$
and such that any torsion element corresponds to a central symmetry of $\HH^3$.
\end{enumerate}
\item\label{thm:main:2:c} The orientation preserving subgroup of the Coxeter group $\Gamma_{120}^+\acts\HH^4$ of the regular right angled 120-cell satisfy this property.
\end{enumerate}

\end{thm}

Theorem \ref{thm:main}(\ref{thm:main:2:b})
and the following corollary was obtained by Kapovich in \cite{K1}.

\begin{thm}{Corollary}\label{cor:M/J-hyperbolic}
Any finitely presented group $G$
is isomorphic to the fundamental group of $M/J$,
where $M$ is either
\begin{enumerate}[(a)]
\item complete noncompact 3-dimensional oriented  hyperbolic manifold, or
\item compact 3-dimensional oriented hyperbolic manifold with convex boundary.
\end{enumerate}
In both cases $J\:M\to M$ is an isometric involution which has only finite number of fixed points.
\end{thm}

Kapovich used this statement in the proof of the following result, which is discussed in Section~\ref{sec:kapovich}.

\begin{thm}{Kapovich's theorem}\label{kapovich} Let $G$ be a finitely-presented group. Then there exists a $2$-dimensional
irreducible complex-projective variety $W$ with the fundamental group $G$, so that the
only singularities of $W$ are normal crossings and Whitney umbrellas.
\end{thm}

\parit{Comments.}
It is well-known that fundamental groups of complex projective  (or compact K\"ahler)
manifolds satisfy many restrictions,
see e.g. \cite{ABCKT}.
On the other hand varieties that are unions
of collections of coordinate planes in complex projective spaces
can have arbitrary fundamental group.
Such varieties are of course
reducible and have very bad singularities.

The theorem above strengthens the result of Simpson in \cite{simpson} which states that
every finitely-presented group $G$ appears as the fundamental group of a singular
{\it irreducible} complex-projective variety.

\medskip

Further,
from Theorem \ref{thm:main}(\ref{thm:main:2:a}), we get the following result;
it was announced by Aitchison,
but he did not wrote the proof.

\begin{thm}{Aitchison's theorem}\label{cor:M/J-compact}
Any finitely presented group $G$
is isomorphic to the fundamental group of $M/J$,
where $M$ is a closed
 oriented 3-dimensional manifold
and $J\:M\to M$ is a smooth involution which has only isolated fixed points.
\end{thm}
The last result might look surprising since there many restrictions on the fundamental group
of 3-manifolds, see for example \cite{wilton}.

\medskip

Theorem~\ref{thm:main}(\ref{thm:main:2:c}) has two applications described in Section~\ref{sec:dim=3C}.
The first one is a new proof of a
result of Taubes in \cite{T}
concerning complex 3-manifolds
and the second is on symplecic Calabi--Yau manifolds
obtained by the first author and Fine in \cite{FP1}.

\begin{thm}{Taubes' theorem}
\label{thm:taubes}
 For every finitely presented group $G$ there
exists a smooth compact complex $3$-manifold $W$ such that $\pi_1 W=G$.
\end{thm}

Again, an analogous theorem does not hold for K\"ahler manifolds (see \cite{ABCKT}).
In fact all the manifolds obtained by our construction are non-K\"ahler;

\begin{thm}{Symplectic Calabi--Yau manifolds}\label{SCY}
For every finitely presented group $G$ there
exists a smooth compact symplectic $6$-manifold $M^6$
with vanishing first Chern class and the fundamental group isomorphic to $G$.
\end{thm}

Such a diversity of symplectic $6$-folds with $c_1=0$ is quite surprising
since only two types of symplectic $4$-manifolds with
$c_1=0$ are known, namely $K3$-surfaces and $T^2$ bundles over $T^2$
(see for example \cite{tian}).

\section{Open problems}\label{sec:open}

The following question is the main motivation for our construction,
it is stated on page 12 in \cite{gromov}.

\begin{thm}{Gromov's question}
Is it true that every compact smooth manifold $M$
is PL-homeomorphic
to the underlying space of a hyperbolic orbifold?

In other words, is there a discrete co-compact isometric action on the hyperbolic space with the quotient space PL-homeomorphic to $M$?
\end{thm}

The question is open if $n\ge 4$.
In dimensions $2$ and $3$ the answer is yes.
Moreover, in these dimensions one can fix a compact hyperbolic orbifold $\mathcal{O}$
and get all the manifolds by passing to finite orbicovers of $\mathcal{O}$.
In the 2-dimensional case one can take the orbifold formed by the hyperbolic triangle with angles $\tfrac\pi2$, $\tfrac\pi3$ and $\tfrac\pi{5{\cdot}6{\cdot}7}$; this is easy to prove.
In the 3-dimensional case one can take right angled hyperbolic dodecahedron;
the later is the baby case of theorems proved by Thurston, Hilden, Lozano, Montesinos and Whitten;
see \cite{HLMW} and the references there in.

\medskip

Let us state more specific open questions.

\begin{thm}{Question}
Is it true that the underlying space of any compact 1000-dimensional hyperbolic orbifold has nontrivial fundamental group?
\end{thm}

Clearly the answer ``yes'' would imply ``no'' for the Gromov's question.
The following is a yet more specific conjecture,
if true, it also gives a negative answer to  Gromov's question.

\begin{thm}{Conjecture}
The underlying space of a 1000-dimensional hyperbolic orbifold
can not be homeomorphic to the 1000-dimensional disc nor to the 1000-dimensional sphere.
\end{thm}

The only approach to this conjecture that we see is to generalize Vinberg's proof of non-existence of hyperbolic cocompact isometric actions on hyperbolic space  of large enough dimension generated by the reflections in the hyperplanes, see \cite{V1} and \cite{V2}.

Let us explain why these two problems are relevant.
The quotient space $\HH^m/\Gamma$ is PL-homeomorphic to a manifold
if and only if the isotropy groups of $\Gamma$ are generated by rotations around subspaces of codimension 2; the later is proved by Mikhailova in \cite{mikhailova}.\footnote{If one exchanges \emph{PL-homeomorphism} to \emph{homeomorphism}
then the answer is more complicated,
in this case the classification of isotropy groups was given recently by Lange in \cite{lange}.}
For manifolds with boundary one needs to add reflections in hyperplanes.

\medskip

An affirmative answer to the following statement was announced by Aitchison,
but later he took it back.
It is closely related to the Aitchison's theorem (\ref{cor:M/J-compact}) and  the Corollary~\ref{cor:M/J-hyperbolic}.
A modification of our construction might lead to a ``yes'' answer;
 we do not see any approach for a ``no'' answer.

\begin{thm}{Question}
Is it true that any finitely presented group appears as the fundamental group of a quotient $M/J$, where $M$ is a compact 3-dimensional hyperbolic manifold and $J\:M\to M$ is an isometric involution?

If the answer is ``yes'' can we assume in addition that $J$ has only isolated fixed points?
\end{thm}

\section{Telescopic actions}\label{sec:telescope}

\begin{thm}{Definition}\label{def:telescopic}
A co-compact properly discontinuous isometric group action $\Gamma\acts X$
on a length-metric space $X$
is called \emph{telescopic} if
given a finitely presented group $G$, there exists
a subgroup $\Gamma'<\Gamma$ of finite index
such that $G$ is isomorphic to the fundamental group of the quotient space $X/\Gamma'$.
\end{thm}

We will construct telescopic actions on particularly nice spaces $X$.
In the base example $X$ will be a 2-dimensional polyhedral $\CAT[-1]$ space;
it will be used to construct telescopic actions on hyperbolic spaces of dimension from $3$ and $4$.

\medskip

A \emph{good compact orbihedron}
is an alternative way to think about
a cocompact group action.
These two languages can be easily translated from one to the other;
for example, passing to a subgroup of finite index
is the same as a finite orbicover of the orbihedron.
It turns out that the orbi-language suites better our purposes.
The definition above can be reformulated in the following way.

\begin{thm}{Definition}\label{def:telescopic-orbi}
A good compact orbihedron $\mathcal{O}$ is called \emph{telescopic}
if given a finitely presented group $G$, there exists
a finite orbicover $\mathcal{O}'\to\mathcal{O}$
such that $G$ is isomorphic to the fundamental group
of the underlying space $|\mathcal{O}'|$ of $\mathcal{O}'$.
\end{thm}

Now we turn to the construction of the 2-dimensional telescopic orbihedron $Y$. First we construct the underlying space of $Y$
and
then we equip it with a metric and an orbi-structure.

\parbf{Underlying space.} Consider figure eight $F_8$ with the loops $r$ and $g$ ($r$ is for ``red'' and $g$ is for ``green'').
Let us attach to $F_8$ four discs $\mathcal{B}$, $\mathcal{W}$, $\mathcal{G}$, $\mathcal{R}$
(named for ``black'', ``white'', ``green'', and ``red'')
along $g{*}r^{-1}$, $g{*}r$, $g$ and $r$ respectively.
Denote by $|Y|$ the obtained topological space.

In other words, $|Y|$ is homeomorphic to $\RP^2$
with two discs attached along two lines.
The lines cut two discs from $\RP^2$ which are colored in black and white
and the attached discs are red and green.

\parbf{Orbi-structure.}
Fix some integer $k\ge 2$.

In the interior of each 2-cell of $Y$,
choose $k$ ``singular'' points;
so in total we have $4\cdot k$ singular points.
Let us assume that each singular point is modeled on the singularity $\RR^2/\ZZ_2$,
and the rest of the points are regular;
i.e., they have trivial isotropy groups.
This defines an orbi-structure on $Y$.

Let us equip $Y$ with a metric.
To do this, prepare a pair of right angled $(k+2)$-gons for each $2$-cell of $Y$.
(The $(k+2)$-gon has to be hyperbolic if $k\ge 3$
and has to be a Euclidean rectangle if $k=2$.)
Gluing the pair of $(k+2)$-gons
along the corresponding $k+1$ sides
we get a disc.
The pair of unglued sides form its boundary of the disc
and the $k$ vertices  which appear in its interior correspond to the singular points.
For a right choice of size of these $(k+2)$-gons the boundaries of four disks can be glued together by length-preserving maps.

The orbihedron $Y$ admits the universal orbi-cover $X$
which is double branching at each singular point.
The induced metric on $X$ is $\CAT[0]$.
Moreover, in the case  $k\ge 3$ we can make it to be $\CAT[-1]$.
Denote by $\Gamma\acts X$ the group of its deck transformations;
this will be the action corresponding to $Y$;
in particular $Y$ is a good orbifold; see \cite{haefliger} for more details.

\parbf{Orbicovers of $\bm{Y}$.}
Consider a two-dimensional CW-complex $W$ which
satisfies the following three conditions.

\begin{enumerate}
\item
The two cells of $W$ can be colored in 4 colors black, white, green and red,
in such a way that black-and-white cells form a connected surface $\Sigma$
(which is not orientable in general).
The red and green cells are attached to $\Sigma$ along a collection of curves which will be called red and green correspondingly.

\item Each curve (green or red) intersects with at least one other curve
 and the intersections are transversal in $\Sigma$.
Two curves of the same color can not intersect.

\item\label{lem:X:discs}
The red and green curves cut $\Sigma$ into black and white discs in the checkerboard order.
Moreover there is an orientation on each curve such that
if one moves around the boundary of white (black) disc
then red and green segments have the same (correspondingly the opposite) orientation.
\end{enumerate}

\begin{wrapfigure}{r}{43mm}
\begin{lpic}[t(-4mm),b(-2mm),r(0mm),l(0mm)]{pics/W(1)}
\lbl[b]{29.5,30,-85;{\small red}}
\lbl[tr]{17,21.5,16;{\small gr}}
\lbl[tl]{17.2,21.5,0;{\small een}}
\lbl[b]{40,12.7,-16;{\small red}}
\lbl[]{40,30;{\large $\Sigma$}}
\end{lpic}
\end{wrapfigure}

In this case, $W$ is homeomorphic to a finite orbicover of $Y$.
Indeed, it is easy to see that the 1-skeleton of $W$ is a cover of $F_8$
which respects the color and the orientation.
This covering can be extended to a ramified covering $W\to Y$
which is branching only at two given
interior points in each cell with order at most 2.
The later follows since any cover $\SS^1\to\SS^1$ can be extended to
a ramified covering $\DD^2\to \DD^2$
which is branching only at the given
two interior points in each cell with the order at most 2; cf. \cite[Proposition 1]{feighn}.

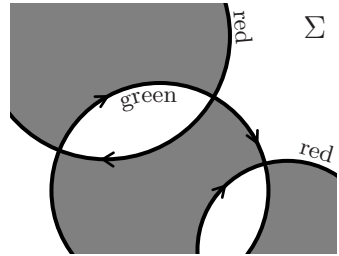
\parbf{Telescopicness.}
It is well known that any finitely presented group $G$
can appear as the fundamental group of a finite two-dimensional CW-complex, say $W$.
This construction enjoys a lot of freedom which can be used to show that $W$ can be chosen so that it satisfies the three conditions above;
see \cite{panov-petrunin} for more details.

Since these conditions imply that $W$ appears us the underlying space of an orbicover of $Y$, we get that $Y$ is telescopic.

\medskip

Note that we proved  \ref{thm:main}(\ref{thm:main:0}).
Indeed, we can assume that $X$ is a $\CAT[-1]$ space.
Therefore given $\gamma\in\Gamma'$,
we have $\gamma\in \Tor \Gamma'$
if and only if $\gamma$ has a fixed point if $X$.
It follows that the fundamental group of $X/\Gamma'$
is isomorphic to the quotient group $\Gamma'/\<\Tor \Gamma'\>$, see \cite{armstrong}.

\parbf{Hyperbolic 3-orbifold.}
Consider a right angled hyperbolic dodecahedron $Q$.
Assume we glue face-to-face a few copies of $Q$.
The obtained space has a natural orbifold structure
if around each vertex of each copy of $Q$
we see the quotient space of $\HH^3$ by a
subgroup of $\ZZ_2\oplus\ZZ_2\oplus\ZZ_2$ generated by reflections in $3$ orthogonal hyperplanes.
The picture is either
$\HH^3$,
$\HH^3/\ZZ_2$,
$\HH^3/(\ZZ_2\oplus\ZZ_2)$
or $\HH^3/(\ZZ_2\oplus\ZZ_2\oplus\ZZ_2)$;
in the second and third cases there exist geometrically different actions.

It turns out that these local rules leave enough freedom to
construct an orbifold $\mathcal{O}$ which underlying space has almost arbitrary big scale structure.
In particular, we can glue an orbifold which looks very much like $Y$.

Moreover, this way we can produce a hyperbolic orbifold $\mathcal{O}$
with an embedding $\iota\:|Y|\to|\mathcal{O}|$
such that the following property holds.
For any orbi-cover $Y'\to Y$
there is an orbi-cover $\mathcal O'\to O$
and an embedding
$\iota'\:|Y'|\hookrightarrow|\mathcal O'|$ such that
the following diagram is commutative
$$\raisebox{-0.9cm}{$\begindc{\commdiag}[10]
\obj(0,3)[aa]{$|Y'|$}
\obj(4,3)[bb]{$|\mathcal O'|$}
\obj(0,0)[cc]{$|Y|$}
\obj(4,0)[dd]{$|\mathcal O|$}
\mor{aa}{bb}{$\iota'$}
\mor{aa}{cc}{}
\mor{bb}{dd}{}
\mor{cc}{dd}{$\iota$}
\enddc$}
$$
and $\iota'$ induces an isomorphism
$\pi_1|Y'|\to\pi_1|\mathcal O'|$.

This property is sufficient to conclude that the $\mathcal O$ is telescopic.

\parbf{Variations of the construction.}
In the above construction the isotropy groups of $\mathcal O$
can be equal to any subgroup of the action generated by the reflections
in the coordinate hyperplanes.
In fact,
the freedom in the construction makes possible to get more control on the isotropy groups.

One can leave only the orientation preserving actions
and the action of $\ZZ_2$ by the central symmetry.
This proves Theorem~\ref{thm:main}(\ref{thm:main:2:a}).
Note that the quotient spaces for the orientation persevering subgroups are topologically manifolds.
The underlying space of the obtained orbifold is a manifold with
singular points which admit neighborhoods homeomorphic to the cone over $\RP^2$;
passing to the orienting double cover, we get Aitchison's theorem~\ref{cor:M/J-compact}.

In the same way one can build an orbifold which
has singularities formed by the action of $\ZZ_2$ by central symmetries while the rest of singularities form the topological boundary of underlying space of the orbifold.
In this case the underlying space is the quotient of a convex set $L$
in $\HH^3$
formed by the union of copies of $Q$
by a cocompact action of a group in which the only torsion elements
are central symmetries.
This proves \ref{thm:main}(\ref{thm:main:2:c}) and implies Corollary \ref{cor:M/J-hyperbolic}.

In the $4$-dimensional case one can use
a hyperbolic right-angled 120-cell the same way as we used $Q$.

Similar construction might work in the $n$-dimensional case
once there is a compact $n$-dimensional hyperbolic orbifold.
No such examples are known for big enough $n$.

\section{Complex and symplectic manifolds}
\label{sec:dim=3C}

In this section we explain how Theorem \ref{thm:main}(\ref{thm:main:2:c}) leads to a simple proof of
Taubes' theorem and to a construction of six-dimensional
symplectic Calabi--Yau manifolds with arbitrarily fundamental group.

Both proofs relay on the twistor construction which we are about to explain.

\parbf{Twistor space.}
Let $X$ be an oriented $4$-dimensional
Riemannian manifold.
Then the twistor space of $X$ is defined as the $\SS^2$-bundle over
$X$ whose fiber over a point $p\in X$ consists of all orthogonal operators $J$ acting
on the tangent space $\T_p X$ such that $J^2=-1$ and $J$ induces on $X$ the correct orientation.

We will denote the twistor space of $\HH^4$ by $Z$. All twistor spaces carry a natural almost complex structure, and in the case of the twistor space of $\HH^4$ this structure
is integrable, see \cite{besse}. 
To get an idea of this complex
structure one can conformally identify $\HH^4$ with a half of the round $\SS^4$.
In this case $Z$ can be seen as a {\it ``half''} of the twistor space of $\SS^4$ which is $\CP^3$.

In order to get a symplectic structure $w$ on $Z$
one can identify $Z$ with a particular six-dimensional co-adjoint orbit of
$\SO(1,4)$, \cite[Section 2.3.3]{FP}. This  orbit consist of matrixes $A$ from
$\mathfrak{so}(1,4)$ whose kernel is generated by a positive vector $v\in \RR^{1,4}$
and such that $A^2|_{v^{\perp}}=-\id$.

Note that these complex and symplectic structures
do not give rise to a K\"ahler form on $Z$ since $w$ is not positive on all
complex directions.

The group of orientation
preserving isometries $\Iso_+(\HH^4)$ of $\HH^4$ lifts to the action on its twistor space $Z$
and the action preserves both complex and symplectic structures.

\parbf{About the proofs.}
Applying Theorem \ref{thm:main}(\ref{thm:main:2:c}),
we get a co-compact subgroup $\Gamma'< \Iso_+(\HH^4)$
such that $\pi_1(\HH^4/\Gamma')\cong G$.
In this case the topological fundamental group of
the quotient of the twistor space
$Z/\Gamma'$ equals $G$ as well,
since the fibers of the projection $Z/\Gamma'\to \HH^4/\Gamma'$ are two-dimensional spheres.
The quotient space $Z/\Gamma'$ is both a complex and symplectic orbifold, moreover with respect to the symplectic
structure $c_1(Z/\Gamma')=0$.

In order to prove Theorems \ref{thm:taubes} and \ref{SCY} one has to
resolve the singularities of $Z/\Gamma'$ with respect to complex and symplectic structures.

The complex case is straightforward, the resolution exists by Hironaka theorem;
it also can be constructed explicitly.
The fundamental group does not change after such a resolution and we obtain
a smooth complex 3-dimensional manifold with the fundamental group $G$.
This proves Taubes' theorem.

The symplectic case is more involved since the techniques of resolutions of symplectic
singularities are developed poorly and moreover
the resolution has to preserve the condition $c_1=0$.

The existence of such a symplectic resolution for $Z/\Gamma'$  is proved in  \cite{FP1}.
The proof uses the fact that all the
singularities of the constructed orbifold
are locally modeled on the quotient of $\mathbb C^3$
with the standard symplectic form by a linear action of $\mathbb Z_2$ or
$\mathbb Z_2\oplus \mathbb Z_2$.
As a consequence  Hamiltonian
actions of tori $\TT^1$ and $\TT^3$ can be defined at neighborhoods of singularities
that are compatible in a specific way.
This permits one to apply methods
of toric geometry and use a version of symplectic cutting to resolve the singularities.

\section{Projective surfaces with controlled singularities}
\label{sec:kapovich}

In this section we give a rough outline the proof
of \ref{kapovich}.

\parit{Step 1; Dirichlet tessellation.}
Starting with $\Gamma$ provided by Corollary~\ref{cor:M/J-hyperbolic} one gets an
infinite tessellation $D_x(\Gamma)$ of $\HH^3$ by Dirichlet domains of the action of
$\Gamma$ on $\HH^3$ with respect to a generic point $x$.
Since the action is convex co-compact all domains have a finite number of faces.

\parit{Step 2; weak simplicity.}
For a generic  $x$ the tessellation $D_x(\Gamma)$ is simple outside of it vertices, see \cite[Theorem 1.6]{K1}.
In other words each edge of the tessellation
belongs to exactly three domains.
This property is called \emph{weak simplicity}.

This is a key moment in the proof.
The proof relies on the fact that the only torsion elements in the group provided by Corollary~\ref{cor:M/J-hyperbolic}
are central symmetries.%
\footnote{It is conjectured that for such actions, the tessellation $D_x(\Gamma)$ is simple for a generic  $x$; see
\cite[Conjecture 1.5]{K1}.}

\parit{Step 3; complexification.}  Next one finds a finite index normal
subgroup $\Gamma'\vartriangleleft \Gamma$
such that the action $\Gamma'\acts \HH^3$ is free.
In this case $D_x(\Gamma)/\Gamma'$ is glued from a finite number of hyperbolic
polytopes all equal to the Dirichlet domain of $\Gamma$.

Further the polyhedral complex $D_x(\Gamma)/\Gamma'$ is complexified.
This operation is quite sophisticated and we  explain just its first approximation.

Each polytope $P$ of the complex $D_x(\Gamma)/\Gamma'$ can be realized
as embedded in $\CP^3$ via a natural chain of embeddings
\[P \hookrightarrow \HH^3\hookrightarrow \RP^3\hookrightarrow \CP^3.\]
The two-faces of $P$ span a collections of $\CP^2$'s in $\CP^3$.
Now one can glue all $\CP^3$'s  along these $\CP^2$'s
by  extending linearly the isometric identification
of polyhedrons faces.

In the actual complexification in \cite{K1} before gluing the copies of $\CP^3$ one has
to blow them up in a carefully chosen way.

The obtained variety $\cal P$ is reducible;
it can be shown to be complex projective
with normal crossing singularities by methods developed in \cite{KK}.
The proof relies on the fact that the complex
$D_x(\Gamma)/\Gamma'$ is simple outside of its vertices.

\parit{Step 4; taking quotient.}  Next one considers the quotient of $\cal P$ by $\Gamma/\Gamma'$,
and shows that its fundamental group equals to $G$.
This quotient is an irreducible projective variety
but it is not anymore a variety with normal crossing
singularities.
Additional singularities come from fixed points of
the action of $\Gamma/\Gamma'$.

\parit{Step 5; reducing dimension to two.}
Finally one takes a
hyperplane section of the obtained projective variety.
It has all the deserted
properties, and Whitney umbrellas appear in it because of the above
additional singularities.


\begin{thebibliography}{References}

\bibitem{ABCKT}
Amor\'os, J.; Burger, M.; Corlette, K.; Kotschick, D.; Toledo, D.
Fundamental groups of compact K\"ahler manifolds. (English summary)
Mathematical Surveys and Monographs, 44. American Mathematical Society, Providence, RI, 1996. xii+140 pp. ISBN: 0-8218-0498-7

\bibitem{armstrong}  Armstrong, M. A. The fundamental group of the orbit space of a discontinuous group, Proc. Cambridge Philos. Soc. 64 (1968) 299--301.

\bibitem{besse} Besse, A.
Einstein manifolds.
Reprint of the 1987 edition.
Classics in Mathematics.
Springer-Verlag,
Berlin, 2008.
xii+516 pp.
ISBN: 978-3-540-74120-6


\bibitem{feighn}  Feighn, M. Branched covers according to J. W. Alexander,
Collect.
Math. 37 (1986), no. 1, 55--60.


\bibitem{FP}  Fine, J.;  Panov, D. Hyperbolic geometry and non-K\"ahler
manifolds with trivial canonical bundle. {\it Geometry and Topology}
14 (2010) 1723--1763.

\bibitem{FP1} Fine, J.; Panov, D. The diversity of symplectic Calabi--Yau six-manifolds.
Preprint, arXiv:1108.5944.

\bibitem{gromov}  Gromov, M.
Manifolds: Where Do We Come From? What Are We? Where Are We Going.
Preprint 2010.

\bibitem{haefliger}Haefliger, A.
Orbi-espaces. Sur les groupes hyperboliques d'apr\`{e}s Mikhael Gromov (Bern, 1988), 203--213,
Progr. Math., 83, Birk\-h\"auser Boston, Boston, MA, 1990. \\
\begin{otherlanguage}{russian}
Перевод: Хэфлигер, А. Орбипространства,
Гиперболические группы по Михаилу Громову,
(М. Мир, 1992), 191--201.
\end{otherlanguage}

\bibitem{HLMW} Hilden, H.;
Lozano, M.;
Montesinos, J.; Whitten, W.
On universal groups and three-manifolds. {\it Invent, math.} 87, 441--456 (1987)

\bibitem{K1}  Kapovich, M.
Dirichlet fundamental domains and complex projective varieties.  Invent. math,
January 2013.

\bibitem{KK} Kapovich, M.; Kollar, J. Fundamental groups of links of isolated singularities, Journal AMS (to appear).


\bibitem{lange} Lange, C.
When is the underlying space of an orbifold a topological manifold?
arXiv:1307.4875

\bibitem{tian}  Li, T.-J. Quaternionic vector bundles and Betti numbers of symplectic 4-manifolds with Kodaira dimension zero. {\it Internat. Math. Res. Notices}, 1--28,  (2006).


\bibitem{panov-petrunin} Panov, D.; Petrunin, A. Telescopic actions. Geom. Funct. Anal. 22 (2012), no. 6, 1814–1831.

\bibitem{mikhailova} \begin{otherlanguage}{russian}
Михайлова, М. А. О факторпространстве по действию конечной группы, порожденной псевдоотражениями
{\it Изв. АН СССР. Сер. матем.}, 48:1 (1984),  104--126.\end{otherlanguage}

\bibitem{simpson} Simpson, C. Local systems on proper algebraic V-manifolds, Pure and Applied Math. Quarterly,
7 (2011), p. 1675--1760.

\bibitem{T}  Taubes, C.
The existence of anti-self-dual conformal structures.
{\it J. Differential Geom.} 36 (1992), no. 1, 163--253.




\bibitem{V1}   \begin{otherlanguage}{russian}Винберг, Э.  Б. Отсутствие кристаллографических групп отражений в пространствах Лобачевского большой размерности. \textit{Труды ММО}, 1984, 47,  68--102. \end{otherlanguage}

\bibitem{V2} \begin{otherlanguage}{russian}Винберг, Э.  Б. Гиперболические группы отражений.
\textit{УМН}, 40:1(241) (1985),  29--66;\end{otherlanguage}\\
Translation: Vinberg, E, B. Hyperbolic reflection groups.
Russ. Math. Surv. 40, No.1, 31--75 (1985);

\bibitem{wilton} Wilton, H. 3-manifold groups are known, right?\\
\href{http://ldtopology.wordpress.com/2010/01/26/3-manifold-groups-are-known-right/}{\texttt{ldtopology.wordpress.com}}

\end{thebibliography}
\end{document}